\documentclass[11pt,reqno]{amsart}
\usepackage{amsmath, amssymb, amsthm}

\renewcommand{\(}{\left\(}
\renewcommand{\)}{\right\)}
\renewcommand{\[}{\left\[}
\renewcommand{\]}{\right\]}

\numberwithin{equation}{section}
 \theoremstyle{plain}
\newtheorem{theorem}{Theorem}[section]
\newtheorem{lemma}[theorem]{Lemma}
\newtheorem{remark}[]{Remark}
\newtheorem{definition}[]{Definition}

\newtheorem{corollary}[theorem]{Corollary}
\newtheorem{proposition}[theorem]{Proposition}
\newtheorem{example}[]{Example}
\newtheorem{problem}[]{Problem}

   \makeatletter
\def\proof{\@ifnextchar[{\@oproof}{\@nproof}}
\def\@oproof[#1][#2]{\trivlist\item[\hskip\labelsep\textit{#2 Proof of\
#1.}~]\ignorespaces}
\def\@nproof{\trivlist\item[\hskip\labelsep\textit{Proof.}~]\ignorespaces}

\makeatother
\begin{document}
\title[A new generalization of the minimal excludant]{A new generalization of the minimal excludant arising from an analogue of Franklin's identity}
\author{Subhash Chand Bhoria}
\address{Subhash Chand Bhoria, Department of Mathematics, Pt. Chiranji Lal Sharma Government College, Urban Estate, Sector-14, Karnal, Haryana - 132001, India.}
\email{scbhoria89@gmail.com}

\author{Pramod Eyyunni}
\address{Pramod Eyyunni, Discipline of Mathematics,
Indian Institute of Technology Indore, Simrol, Indore, Madhya Pradesh - 453552, India.} 
\email{pramodeyy@iiti.ac.in, pramodeyy@gmail.com}

\author{Bibekananda Maji}
\address{Bibekananda Maji, Discipline of Mathematics,
Indian Institute of Technology Indore, Simrol, Indore, Madhya Pradesh - 453552, India.}
\email{bibekanandamaji@iiti.ac.in}

\thanks{$2020$ \textit{Mathematics Subject Classification.} Primary 11P81, 11P84; Secondary 05A17. }

\begin{abstract}
Euler's classical identity states that the number of partitions of an integer into odd parts and distinct parts are equinumerous. Franklin gave a generalization by considering partitions with exactly $j$ different multiples of $r$, for a positive integer $r$. We prove an analogue of Franklin's identity by studying the number of partitions with $j$ multiples of $r$ in total and in the process, discover a natural generalization of the minimal excludant (mex) which we call the $r$-chain mex. Further, we derive the generating function for $\sigma_{rc} \textup{mex}(n)$, the sum of $r$-chain mex taken over all partitions of $n$, thereby deducing a combinatorial identity for $\sigma_{rc} \textup{mex}(n)$, which neatly generalizes the result of Andrews and Newman for $\sigma \textup{mex}(n)$, the sum of mex over all partitions of $n$.    
\end{abstract}  
\maketitle
\section{Introduction}
For a natural number $n$, a partition is a non-increasing sequence of positive integers that sum to $n$. The function $p(n)$ counts the number of partitions of $n$. For instance $p(3) = 3$, corresponding to the partitions $(3), \ (2, 1), \ (1, 1, 1)$ of $3$. Partitions are also represented as a finite sum or by noting the `frequencies' of the integers in the partition. Thus, the partitions of $3$ can also be written as $3, \ 2+1, \ 1+1+1$ or $\left< 3^1\right>, \left< 1^1 \ 2^1\right>, \left< 1^3\right>$, where $a^m$ implies that the integer $a$ repeats $m$ times in the partition. In this work, we usually use the finite sum notation but also invoke the frequency notation as needed.
 
In the mid eighteenth century, Euler proved a classical result linking two types of partition functions. He showed that
\begin{theorem} \label{Euler class}
The number of partitions of a positive integer $n$ into odd parts equals the number of those into distinct parts.
\end{theorem}
For example, the partitions of $6$ into odd parts are $5+1, \ 3+3, \ 3+1+1+1, \ 1+1+1+1+1+1$, and into distinct parts are $6, \ 5+1, \ 4+2, \ 3+2+1$, so there are four partitions of each kind. This set in motion the field of partition identities, with many tantalizing identities being discovered subsequently. For instance, one of the famous pair of Rogers-Ramanujan identities states that the number of partitions of $n$ into distinct parts where the gap between successive parts is at least two equals the number of partitions of $n$ into parts congruent to $1$ or $-1$ modulo $5$. This can be seen as a `modulo 5' counterpart to Theorem \ref{Euler class} if we rephrase one half of the latter identity as `partitions into parts congruent to $1$ or $-1$ modulo $4$'.

There have been many generalizations and refinements of Euler's identity. Sylvester, in his magnum opus \cite{sylvester}, proved that the number of partitions of $n$ with $k$ different odd parts occurring equals the number of partitions of $n$ into distinct parts with exactly $k$ (maximal) blocks of consecutive integers. On the other hand, Fine \cite[p. 46]{fine Euler ref} showed that the number of partitions of $n$ into distinct parts with largest part $k$ equals the number of partitions of $n$ into odd parts such that the largest part plus twice the number of parts is $2k+1$. See also Alladi \cite{alladi refine} and the references therein for further details.

A well-known generalization of Euler's identity is an identity by Glaisher, which is as follows:
\begin{theorem} \label{glaisher}
The number of partitions of $n$ into parts not divisible by a positive integer $r$ equals the number of partitions of $n$ with each integer occurring less than $r$ times.
\end{theorem}
Clearly, the case $r=2$ corresponds to Euler's identity, Theorem \ref{Euler class}. Franklin \cite{franklin} provided a remarkable generalization of Glaisher's identity which is as follows:
\begin{theorem}\label{franklin}
The number of partitions of $n$ in which exactly $j$ different parts (these parts can be repeated) are divisible by $r$ equals the number of partitions of $n$ in which exactly $j$ different parts occur at least $r$ times and rest of the parts appear at most $r-1$ times .
\end{theorem}
We now discuss a partition statistic of recent prominence, namely, the minimal excludant which provides an additional insight into our analogue of Franklin's identity (See Theorem \ref{3-way F-analogue} in Section \ref{motiv} below). The minimal excludant (mex) of an integer partition is the least positive integer that does not occur as a part in that partition. This was first introduced by Grabner and Knopfmacher \cite{grabner} under the name of `least gap' of a partition along with an analytic study of its properties. Recently, Andrews and Newman \cite{andrews-newmanI} revived the investigation of the least gap, especially from a combinatorial viewpoint, under the name of minimal excludant. They defined the sum of minimal excludants function, namely, $\sigma \textup{mex}(n) := \sum_{\pi \in \mathcal{P}(n)} \textup{mex} (\pi)$, where $\textup{mex}(\pi)$ denotes the mex of the partition $\pi$ and showed that 
\begin{equation}\label{sigma mex comb id}
\sigma \textup{mex}(n) = D_2(n).
\end{equation} 
Here, $D_2(n)$ represents the number of partitions of $n$ into distinct parts with two colors. They were able to prove this by deriving the generating function of $\sigma \textup{mex}(n)$:
\begin{equation} \label{sigma mex gfn}
\sum_{n=0}^{\infty} \sigma \textup{mex} (n) q^n = (-q; q)_{\infty}^2.
\end{equation}
Here, and throughout the paper, $|q|<1$ and if $a \in \mathbb{C}$, we define
\begin{align*}
(a;q)_0 :&=1, \qquad \\
(a;q)_n :&= (1-a)(1-aq)\cdots(1-aq^{n-1}),
\qquad n \geq 1, \\
(a;q)_{\infty}:&= \lim_{n \to \infty}(a;q)_n.
\end{align*}
In a subsequent paper, Andrews and Newman \cite{andrews-newmanII} defined a generalized mex called $\textup{mex}_{A, a}(\pi)$, which is the smallest positive integer congruent to $a$ modulo $A$ that does not appear in $\pi$. In particular, $\textup{mex}_{1, 1}(\pi)$ coincides with $\textup{mex}(\pi)$. They defined the partition function $p_{A, a}(n)$ to be the number of partitions $\pi$ of $n$ with $\textup{mex}_{A, a}(\pi) \equiv a \pmod{2A}$ and studied its generating function $F_{A, a}(q)$ for certain relations between $A$ and $a$. As a result they deduced several fascinating identities involving the $p_{A, a}(n)$'s, one of them being
\begin{theorem} \label{odd mex non-neg crank}
The number of partitions of $n$ whose mex is odd, that is, $p_{1, 1}(n)$ equals the number of partitions of $n$ with a non-negative crank. 
\end{theorem}
Sellers and da Silva \cite{sellers-da silva} provided parity characterizations for $p_{1, 1}(n)$ and $p_{3, 3}(n)$. Dhar, Mukhopadhyay and Sarma \cite{dhar et al} generalized Theorem \ref{odd mex non-neg crank} to obtain identities between various $p_{A, a}(n)$ functions and partition functions with arbitrary upper and lower bounds on the statistics rank and crank. 

There have been other types of generalizations of mex as well. Hopkins, Sellers and Stanton \cite{HSS} introduced the mex function $\textup{mex}_j(\pi)$ to be the least integer greater than $j$ that is not a part of $\pi$, where $j$ is a positive integer appearing as a part in $\pi$. In the case when $j$ does not occur in $\pi$, $\textup{mex}_j(\pi)$ is left undefined. If $j$ is considered to be a part of every partition, then the usual mex corresponds to the case $j=0$. Then they showed that for $j \geq 0$, the number of partitions of $n$ whose crank $\geq j$ equals the number of partitions $\pi$ of $n$ with an odd value of $\textup{mex}_j (\pi) - j$. Setting $j=0$ gives back Theorem \ref{odd mex non-neg crank}. Ballantine and Merca \cite{least r gap} examined the combinatorial aspects of the least $r$-gap, which is the least integer occurring less than $r$ times in a partition, and also studied the corresponding $\sigma \textup{mex}$ function. For more such generalizations, the reader can see \cite{HSYee, konan}. 

The corresponding sum of excludants functions have been studied for various types of excludants. Andrews and Newman, apart from initiating the study of $\sigma\textup{mex}(n)$, also examined the quantity $\sigma\textup{moex}(n):= \sum_{\pi \in \mathcal{P}(n)} \textup{moex}(\pi)$ with $\textup{moex}(\pi)$ being the smallest odd natural number missing from the partition $\pi$. The present authors, along with Kaur, considered a restricted version of $\sigma\textup{mex}(n)$, namely the function $\sigma_d\textup{mex}(n)$ which adds up the excludants over all partitions of $n$ into distinct parts. They showed that its generating function is related to Ramanujan's function $\sigma(q)$ (see \cite{mex distinct} for more details). 

In the rest of the paper, we use the following notations: (here $\pi$ is an integer partition)
\begin{align*}
\#(\pi):&= \textup{the number of parts in} \ \pi, \\
s(\pi):&= \textup{the smallest part in} \ \pi, \\
\ell(\pi):&= \textup{the largest part in} \ \pi, \\
\nu(r):&= \textup{the frequency of the integer $r$ in} \ \pi, \\
\mathcal{P}(n)&:= \textup{the collection of partitions of} \ n.
\end{align*}

In the next section, we state our main results along with their motivation. Included in the same section is a new generalization of mex, namely $r$-chain mex (see Definition \ref{r-chain mex} in Section \ref{motiv}) and results pertaining to it in Section \ref{r-chain mex subsec}. After that, in Section \ref{prelim}, we list some well-known results used in the sequel. Then in Section \ref{proofs}, we prove all the main results and finally give pointers to possible directions of research in the last section. 

\section{Motivation and main results} \label{motiv}
The first result that we establish in this paper is an identity relating the mex and certain parts in the partition. It reads as follows:

\begin{proposition}\label{r Euler mexversion}
Let $\alpha(n, j)$ be the number of partitions of $n$ which have exactly $j$ parts greater than their mex and let $e(n, j)$ denote the number of partitions of $n$ which have exactly $j$ even parts. Then we have
\begin{align*}
\sum_{n=0}^{\infty}\sum_{j=0}^{\infty} \alpha(n, j) w^j q^n=\frac{1}{(q;q^2)_{\infty}(wq^2;q^2)_{\infty}} = \sum_{n=0}^{\infty}\sum_{j=0}^{\infty} e(n, j) w^j q^n,
\end{align*}
so that
\begin{equation*} 
\alpha(n, j)=e(n, j).
\end{equation*}
\end{proposition}
From this proposition, it is not hard to see that the total number of even parts in all partitions of $n$ equals the total number of parts greater than mex over all partitons of $n$. Andrews and Merca \cite{universal even} studied the total number of even parts in all partitions into distinct parts and gave new and interesting combinatorial interpretations for it. 
\begin{remark}
The case $j=0$ of Proposition \ref{r Euler mexversion} says that the number of partitions of $n$ with zero parts greater than mex equals the number of partitions into zero even parts. But partitions of $n$ with no parts greater than mex are precisely those in which the smallest part is $1$, with every integer between the smallest and largest parts occuring as a part. We denote this class of partitions of $n$ by  $\mathcal{P}^* (n)$. It is easily seen by conjugation (or) well-known that this set is equinumerous with the set of partitions of $n$ into distinct parts. Thus, the $j=0$ case corresponds to a restatement of Euler's identity.
\end{remark}
Motivated by the above observation, we apply conjugation on partitions with $j$ parts greater than mex to get a partition whose \textit{largest repeating part} is $j$. Thus, we get a new generalization of Euler's result which keeps track of the number of even parts in a partition.
\begin{proposition} \label{reven-Euler}
Let $j$ be a non-negative integer. Then the number of partitions of $n$ with exactly $j$ even parts equals the number of partitions of $n$ whose largest repeating part is $j$, which in turn is the same as the number of partitions of $n$ with $j$ parts greater than mex.
\end{proposition}
For $j=0$, the equality of the first two phrases in Proposition \ref{reven-Euler} is precisely Euler's identity. So Euler's identity, seen in conjunction with this family of results, is actually more about counting even parts in a partition rather than about odd parts, as may appear on a first look. It is interesting to compare Proposition \ref{reven-Euler} with another classical identity, namely, the number of partitions of $n$ into $j$ parts equals the number of partitions of $n$ whose largest part is $j$.
We illustrate Proposition \ref{reven-Euler} below in the case of $n=7$ and $j=2$.

\begin{center}
\begin{tabular}{|c|c|c|}
\hline
Partitions with: $2$ parts $> \textup{mex}$ & $2$ even parts & largest repeating part $2$\\
\hline
$5+2$ & $4+2+1$ & $3+2+2$ \\
\hline
$4+3$ & $2+2+1+1+1$ & $2+2+2+1$ \\
\hline
$3+3+1$ & $3+2+2$ & $2+2+1+1+1$ \\
\hline
\end{tabular} 
\end{center}

Inspired by Franklin's identity in Theorem \ref{franklin}, where partitions with exactly $j$ different parts (these may be repeated) divisible by $r$ were studied, we give a partial generalization of Proposition \ref{reven-Euler} by examining partitions which have a total of $j$ multiples of $r$. 
\begin{proposition}\label{r-chain mex part1}
For a positive integer $r$, define an $r$-\emph{repeating part} of a partition to be an integer which occurs at least $r$ times in that partition. Then the number of partitions of $n$ with exactly $j$ parts divisible by $r$ equals the number of partitions of $n$ in which the largest $r$-repeating part is $j$.
\end{proposition}
Incidentally, Proposition \ref{r-chain mex part1} with $r=1$ is nothing but the fundamental partition identity relating the number of parts and the largest part. Recall that the second equality in Proposition \ref{reven-Euler} was discovered by conjugation of the partitions which had $j$ parts greater than mex. Since conjugation is an involution, a partition with largest repeating part $j$ gives rise to a partition with $j$ parts greater than the mex. This motivates us to apply conjugation to partitions with largest $r$-repeating part $j$ in Proposition \ref{r-chain mex part1}. In the process, we discover a new generalization of mex defined below:
\begin{definition}[$r$-chain mex] \label{r-chain mex}
Let $r$ be a positive integer. The $r$-chain mex of a partition is the least positive integer $k$ such that the integers $k, \ k+1, \ \dots, \ k+r-1$ do not occur as parts in the partition.
\end{definition}
Here, $k+r$ may or may not occur as a part. Clearly, $1$-chain mex is the usual mex of a partition.
\begin{example}
Consider the partition $7+4+4+4+3+1+1$; the $\textup{mex}$ is $2$, the $\textup{2-chain mex}$ is $5$ and for $r \geq 3$, the $\textup{r-chain mex}$ is $8$.
\end{example}
Our first main result is a generalization of Proposition \ref{reven-Euler}:
\begin{theorem} \label{3-way F-analogue}
Suppose that $j$ is a non-negative integer and $r \geq 2$ is a positive integer. Then the number of partitions of $n$ with exactly $j$ multiples of $r$ equals the number of partitions of $n$ whose largest $r$-repeating part is $j$, which is also equal to the number of partitions of $n$ with $j$ parts greater than the $(r-1)$-chain mex.
\end{theorem}
An illustration of Theorem \ref{3-way F-analogue} for $n=7, j=1$ and $r=3$:
\begin{center}
\begin{tabular}{|c|c|c|}
\hline
Partitions with: $1$ part $> 2$-chain mex & $1$ multiple of $3$ & largest $3$-repeating part $1$\\
\hline
$7$ & $6+1$ & $4+1+1+1$ \\
\hline
$6+1$ & $4+3$ & $3+1+1+1+1$ \\
\hline
$5+2$ & $3+2+2$ & $2+2+1+1+1$ \\
\hline
$5+1+1$ & $3+2+1+1$ & $2+1+1+1+1+1$ \\
\hline
$4+1+1+1$ & $3+1+1+1+1$ & $1+1+1+1+1+1+1$ \\
\hline
\end{tabular} 
\end{center}

\begin{remark}
The reader may have observed that though Proposition \ref{r Euler mexversion} follows from Theorem \ref{3-way F-analogue}, we have supplied a proof for it nevertheless. Indeed, its proof is a direct one, whereas in the case of Theorem \ref{3-way F-analogue}, we separately show that both the analogous quantities to those in Proposition \ref{r Euler mexversion} are equal to the number of partitions of $n$ with largest $r$-repeating part $j$.
\end{remark}
\subsection{The $r$-chain mex of an integer partition} \label{r-chain mex subsec}
We proceed to investigate the properties of the $r$-chain mex. In particular, we study the sum of $r$-chain mex taken over all partitions of a postitive integer $n$ and find its generating function along with an elegant combinatorial identity, thereby generalizing the results of Andrews and Newman in \eqref{sigma mex comb id} and \eqref{sigma mex gfn}.
\begin{definition} \label{sigma chain mex}
Let $\textup{rc-mex}(\pi)$ denote the $r$-chain mex of the partition $\pi$. Then the sum of $r$-chain mex taken over all the partitions of $n$ is denoted by $\sigma_{rc} \textup{mex}(n)$, that is,
\begin{equation*}
\sigma_{rc}\textup{mex}(n) :=\sum_{\pi \in \mathcal{P}(n)} \textup{rc-mex} (\pi),
\end{equation*}
where $\mathcal{P}(n)$ is the collection of all partitions of $n$.
\end{definition}
Since the $1$-chain mex of a partition $\pi$ is its usual mex, we have $\textup{1c-mex}(\pi) = \textup{mex}(\pi)$ and hence $\sigma_{1c}\textup{mex}(n) = \sigma\textup{mex}(n)$.

The generating function for $\sigma_{rc}\textup{mex}(n)$ is not derivable in the manner as that of $\sigma \textup{mex}(n)$. The reason being as follows: Suppose $t$ is the $r$-chain mex of a partition. Then $t, t+1, \dots, t+r-1$ are missing from the partition and there is no restriction on parts $\geq t+r$. Note that $t-1$ is necessarily a part (otherwise $t$ will not be the $r$-chain mex), and integers $ \leq t-1$ have to appear in such a way such that the gap between successive parts is at most $r$. This is the issue as the generating function for partitions with the gap between successive parts being at most $r$ (for $r > 1$) is not straightforward. Hence shadowing the approach of Andrews and Newman \cite{andrews-newmanI} for $\sigma \textup{mex}(n)$ is untenable. 

So to understand the sum of $r$-chain mex over all partitions of $n$, we split up the partitons of $n$ according to the number of parts (say $j$) greater than the $r$-chain mex and compute the sum of the corresponding $r$-chain mex and then finally sum over all $j \geq 0$. Before going to the generating function of $\sigma_{rc} \textup{mex}(n)$, we require an auxiliary lemma. By keeping track of what happens to the mex in the conjugation process while proving Theorem \ref{3-way F-analogue}, we have
\begin{lemma} \label{refine lemma}
The number of partitions of $n$ with $j$ parts greater than the $r$-chain mex $k$ equals the number of partitions of $n$ whose largest $(r+1)$-repeating part is $j$, with $k-1$ parts being greater than $j$.
\end{lemma}
We are ready to state the generating function for $ \sigma_{rc}\textup{mex}(n)$:
\begin{theorem}\label{gfn for sigma rc mex}
If $r$ is a postitive integer, then
\begin{equation*} 
\sum_{n=0}^{\infty} \sigma_{rc} \textup{mex} (n) q^n = - \frac{(r-1)}{(q; q)_{\infty}}+ \frac{(q^{r+1}; q^{r+1})_{\infty}}{(q; q)_{\infty}} \sum_{m=1}^{r} \frac{1}{(q^m; q^{r+1})_{\infty}}. 
\end{equation*}
\end{theorem}
This generating function gives us an elegant combinatorial interpretation for $\sigma_{rc} \textup{mex} (n)$ as follows:
\begin{corollary}\label{ID FOR SIGMA RC MEX}
Let $n \geq 0$ and $r \geq 1$. Denote by $p_m(j, n)$, the number of $m$-regular partitions of $n$ in which parts congruent to $j$ modulo $m$ come in two colors. Then
\begin{equation}\label{id for sigma rc mex}
\sigma_{rc} \textup{mex} (n) = -(r-1) p(n) + \sum_{j=1}^{r} p_{r+1} (j, n).
\end{equation}
\end{corollary}
Setting $r=1$ in Theorem \ref{gfn for sigma rc mex} gives us
\begin{align*}
\sum_{n=0}^{\infty} \sigma_{1c} \textup{mex} (n) q^n = \frac{(q^{2}; q^{2})_{\infty}}{(q; q)_{\infty}} \times \frac{1}{(q; q^2)_{\infty}} = \frac{1}{(q; q^2)_{\infty}^2} = (-q; q)_{\infty}^2,
\end{align*}
where the last equality follows from Euler's identity; and this is precisely the generating function for $\sigma \textup{mex}(n)$ mentioned in \eqref{sigma mex gfn}, and proved by Andrews and Newman \cite[Theorem $1.1$]{andrews-newmanI}. Similarly, putting $r=1$ in \eqref{id for sigma rc mex} yields
\begin{equation*}
\sigma_{1c} \textup{mex} (n) = p_2(1, n),
\end{equation*}
the number of $2$-regular partitions of $n$ where parts congruent to $1$ modulo $2$ appear in two colors, or in other words, the number of partitions of $n$ into odd parts with two colors. But, by Euler's identity, this is the same as the number of partitions of $n$ into distinct parts with two colors, $D_2(n)$ and hence we recover \eqref{sigma mex comb id}.

\section{Preliminaries} \label{prelim}
The $q$-binomial theorem is the following useful transformation of an infinite series: (\cite[p.~17, Equation (2.2.1)]{gea1998}) 
\begin{equation}\label{q-binomial}
\sum_{n=0}^{\infty}\frac{(a; q)_n}{(q; q)_n} z^n=\frac{(az; q)_\infty}{(z; q)_\infty}, \quad \text{for} \quad |z| <1.
\end{equation}
A transformation formula of $q$-hypergeometric series due to Heine is given by \cite[p.~359, (III.1)]{gasper}
\begin{align}\label{heine}
{}_{2}\phi_{1}\left( \begin{matrix} a, & b \\
& c \end{matrix} \, ; q, z  \right) = \frac{(b; q)_{\infty}(a z; q)_{\infty }}{(c; q)_{\infty}(z; q)_{\infty }} {}_{2}\phi_{1}\left( \begin{matrix} \frac{c}{b} , & z \\
& a z \end{matrix} \, ; q, b  \right).
\end{align}
We also note down an equivalent formulation for partition conjugation for easy reference in the sequel: (see \cite[pp. 6-7]{sills})
\begin{definition}\label{conj alg}
The conjugate $\lambda^{'} = (\lambda_1^{'}, \lambda_2^{'}, \dots, \lambda_{\lambda_1}^{'})$ of the partition $\lambda = (\lambda_1, \lambda_2, \dots, \lambda_{\#(\lambda)})$ is defined to be the partition in which the $j$th part is given by
\begin{equation*}
\lambda_j^{'} = \sum_{t=j}^{\lambda_1} f_t(\lambda)
\end{equation*}
for $j = 1, 2, \dots, \lambda_1$ and where $f_t(\lambda)$ is the number of times the integer $t$ occurs in $\lambda$. Equivalently, in the `frequency' notation, we can write
\begin{equation}\label{freq form conj}
\lambda^{'} = \left<1^{\lambda_1 - \lambda_2} 2^{\lambda_2 - \lambda_3} 3^{\lambda_3 - \lambda_4} 4^{{\lambda_4 - \lambda_5}} \cdots \right>.
\end{equation}
\end{definition}

\section{Proofs of main results} \label{proofs}
\begin{proof}[Proposition \ref{r Euler mexversion}][] We start with the generating function for $\alpha(n, j)$:
\begin{align}
\sum_{n=0}^{\infty}\sum_{j=0}^{\infty} \alpha(n, j) w^j q^n&=\sum_{m=1}^{\infty}\frac{q^1}{1-q^1}\times\frac{q^2}{1-q^2}\times\cdots\times\frac{q^{m-1}}{1-q^{m-1}} \times \prod_{t=j+1}^{\infty} \frac{1}{1-wq^{t}}\nonumber
\\&=\sum_{m=1}^{\infty}\frac{q^{m\choose2}}{(q;q)_{m-1}(wq^{m+1};q)_{\infty}}\nonumber
\\&=\frac{1}{(wq^2;q)_{\infty}}\sum_{m=1}^{\infty}\frac{(wq^2;q)_{m-1}q^{m\choose2}}{(q;q)_{m-1}}\nonumber
\\&=\frac{1}{(wq^2;q)_{\infty}}\sum_{m=0}^{\infty}\frac{(wq^2;q)_mq^{m+1\choose2}}{(q;q)_m}\nonumber
\\&=\frac{1}{(wq^2;q)_{\infty}}\lim_{A\rightarrow0}\sum_{m=0}^{\infty}\frac{(-q/A;q)_m(wq^2;q)_m}{(0;q)_m(q;q)_m} A^m \nonumber
\\&=\frac{1}{(wq^2;q)_{\infty}}\lim_{A\rightarrow0}{}_{2}\phi_{1}\left( \begin{matrix} -q/A, & wq^2 \\
& 0 \end{matrix} \, ; q, A  \right) \label{Lim_A_in RHS}.
\end{align}
Using \eqref{heine} with $a = -q/A, \ b=wq^2, \ c=0$ and $z=A$ in \eqref{Lim_A_in RHS}, we get
\begin{align}
\sum_{n=0}^{\infty}\sum_{j=0}^{\infty} \alpha(n, j) w^j q^n &=\frac{1}{(wq^2;q)_{\infty}}\lim_{A\rightarrow0}{}_{2}\phi_{1}\left( \begin{matrix} -q/A, & wq^2 \\
& 0 \end{matrix} \, ; q, A  \right)\nonumber
\\&=\frac{1}{(wq^2;q)_{\infty}}\times \lim_{A\rightarrow 0}\frac{(wq^2;q)_{\infty}(-q;q)_{\infty}}{(0;q)_{\infty}(A;q)_{\infty}}{}_{2}\phi_{1}\left( \begin{matrix} 0, & A \\
& -q \end{matrix} \, ; q, wq^2  \right)\nonumber
\\&=(-q;q)_{\infty} {}_{2}\phi_{1}\left( \begin{matrix} 0, & 0 \\
& -q \end{matrix} \, ; q, wq^2  \right)\nonumber
\\&=(-q;q)_{\infty} \sum_{m=0}^{\infty}\frac{w^mq^{2m}}{(-q;q)_m(q;q)_m}\nonumber
\\&=(-q;q)_{\infty} \sum_{m=0}^{\infty}\frac{(wq^2)^m}{(q^2;q^2)_m}\nonumber
\\&=\frac{1}{(q;q^2)_{\infty}(wq^2;q^2)_{\infty}} \nonumber,
\end{align}
where the last equality follows from Euler's identity and an application of \eqref{q-binomial} to base $q^2$ with $a=0$ and $z=wq^2$. This completes the proof of the proposition.
\end{proof}

\begin{proof}[Proposition \ref{reven-Euler}][]
It follows from Proposition \ref{r Euler mexversion} and the special case $r=2$ of Theorem \ref{3-way F-analogue}. (see below)
\end{proof}

\begin{proof}[Theorem \ref{3-way F-analogue}][]
We proceed with the proof in two parts: \\
\textbf{Part 1:} We show firstly that the number of partitions of $n$ with $j$ multiples of $r$ equals the number of partitions of $n$ in which the largest $r$-repeating part is $j$. The generating function for the number of partitions with $w$ keeping track of the number of multiples of $r$ in the partition is given by
\begin{equation} \label{multiples of r I}
\frac{(q^r; q^r)_{\infty}}{(q; q)_{\infty}} \times \frac{1}{(wq^r; q^r)_{\infty}} = \frac{(q^r; q^r)_{\infty}}{(q; q)_{\infty}} \sum_{n=0}^{\infty} \frac{(wq^r)^n}{(q^r; q^r)_n},
\end{equation}
the right side resulting from an application of \eqref{q-binomial} with $q$ replaced by $q^r$ and then setting $a=0$ and $z=wq^r$.
Since we are looking for the generating function of the number of partitions with $j$ multiples of $r$, we extract the coefficient of $w^j$ in \eqref{multiples of r I} to get
\begin{align*}
\frac{(q^r; q^r)_{\infty}}{(q; q)_{\infty}} \times \frac{q^{rj}}{(q^r; q^r)_j} &= \frac{q^{rj}}{(q; q)_{\infty}} \left(q^{r(j+1)}; q^r \right)_{\infty} \nonumber \\
&= \frac{q^{rj}}{(q; q)_j} \prod_{m=j+1}^{\infty} \frac{1-q^{mr}}{1-q^m} \nonumber \\
&= \frac{q^{rj}}{(q; q)_j} \prod_{m=j+1}^{\infty} \left( 1 + q^m + \cdots + q^{m(r-1)}\right),
\end{align*}
which is nothing but the generating function for the number of partitions with largest $r$-repeating part $j$ and hence the equality between the two kinds of partitions follows for each postitive integer $n$.

\textbf{Part 2:} Next, we prove that the number of partitions of $n$ with $j$ parts greater than $(r-1)$-chain mex again equals the number of partitions of $n$ with largest $r$-repeating part $j$. We start by considering the generating function for the number of partitions while simultaneously tracking the number of parts greater than $(r-1)$-chain mex. Suppose, for $k \geq 0$, the $(r-1)$-chain mex of a partition $\mu$ is $k+1$. That means the integers $k+1, \ k+2, \ \dots, \ k+r-1$ are missing from $\mu$. Note that $k$ should occur as a part in $\mu$, for otherwise we would have the $(r-1)$-chain mex of $\mu$ to be $k$ or less, a contradiction. Moreover, the positive integers from $1$ to $k$ occur in $\mu$ in such a way that the gap between successive parts is less than $r$, so as not to contradict the fact that $k+1$ is the $(r-1)$-chain mex of $\mu$. 

Therefore, if $F_k^{<r}(q)$ denotes the generating function for partitions with largest part $k$ and gaps between successive parts being less than $r$, then the generating function for the number of partitions whose $(r-1)$-chain mex equals $k+1$, and where $w$ keeps track of the number of parts greater than  $(r-1)$-chain mex is
\begin{equation*} 
F_k^{<r}(q) \times \frac{1}{(wq^{k+r}; q)_{\infty}}.
\end{equation*}   
Summing over all $k \geq 0$, we obtain the generating function for the number of partitions where $w$ keeps tab on the number of parts greater than $(r-1)$-chain mex:
\begin{equation} \label{j parts > rc-mex I}
\sum_{k=0}^{\infty} F_k^{<r}(q) \times \frac{1}{(wq^{k+r}; q)_{\infty}} = \sum_{k=0}^{\infty} F_k^{<r}(q) \sum_{n=0}^{\infty} \frac{(wq^{k+r})^n}{(q; q)_n},
\end{equation}
where we have invoked \eqref{q-binomial} with $a=0, \ z=wq^{k+r}$.
As we want the generating function of the number of partitions with $j$ parts greater than $(r-1)$-chain mex, we pick the coefficient of $w^j$ in \eqref{j parts > rc-mex I} and get
\begin{equation}\label{j parts > rc-mex II}
\sum_{k=0}^{\infty} F_k^{<r}(q) \frac{(q^{k+r})^j}{(q; q)_j} = \frac{q^{rj}}{(q; q)_j} \sum_{k=0}^{\infty} F_k^{<r}(q) q^{kj}.
\end{equation}
Observe that as $k$ is the largest part in the partitions enumerated by $F_k^{<r}(q)$, the sum on the right side above is the generating function for partitions with largest part repeating more than $j$ times and gaps between successive parts is less than $r$. If $\pi = (\pi_1, \ \pi_2, \ \dots, \ \pi_t)$ is such a partition, then $\nu(\ell(\pi)) > j, \ \pi_i - \pi_{i+1} < r$ for $i = 1, \ \dots, \ t-1$ and $\pi_t < r$. 
Therefore, the conjugate $\pi'$ of $\pi$ satisfies the following properties: $s(\pi') > j$ and by using \eqref{freq form conj}, that is,
\begin{equation*}
\pi' = \left< 1^{\pi_1 - \pi_2} 2^{\pi_2 - \pi_3} \cdots \right>,
\end{equation*}
the gap conditions on $\pi$ give us that each integer in $\pi'$ occurs less than $r$ times. Since conjugation is a bijection, the sum on the right side of \eqref{j parts > rc-mex II} is also equal to the generating function for partitions $\lambda$ with $s(\lambda) > j$ and each integer occurring less than $r$ times, and this is given by $\prod_{t=j+1}^{\infty} \left( 1+q^t + \cdots + q^{t(r-1)}\right)$. Substituting this into \eqref{j parts > rc-mex II} for the sum on the right side there, we obtain the generating function of the number of partitions with $j$ parts greater than $(r-1)$-chain mex as
\begin{equation*}
\frac{q^{rj}}{(q; q)_j}\prod_{t=j+1}^{\infty} \left( 1+q^t + \cdots + q^{t(r-1)}\right),
\end{equation*}
and this being the generating function for the number of partitions with largest $r$-repeating part $j$, the proof of the theorem is complete.
\end{proof}
\begin{proof}[Lemma \ref{refine lemma}][]
We do some bookkeeping with regards to the mex in the proof of Theorem \ref{3-way F-analogue}. By the notation used in the proof of Theorem \ref{3-way F-analogue}, the generating function for the number of partitions with $j$ parts greater than the $r$-chain mex $k$ is the coefficient of $w^j$ in $F_{k-1}^{< (r+1)}(q)  \times \displaystyle\frac{1}{(wq^{k+r}; q)_{\infty}}$, that is
\begin{equation} \label{lemma inter}
F_{k-1}^{< (r+1)}(q) \times \frac{q^{(k+r)j}}{(q; q)_{j}},
\end{equation}
after an appeal to \eqref{q-binomial}. Now, \eqref{lemma inter} can be rewritten as
\begin{equation} \label{lemma inter 2}
\frac{q^{(r+1)j}}{(q; q)_j} \times \left( q^{(k-1)j} \,  F_{k-1}^{< (r+1)}(q) \right) = A(q) \times B(q) \ (\textup{say}).
\end{equation}
In \eqref{lemma inter 2}, $A(q)$ generates partitions $\lambda$ with $\ell(\lambda) = j$ and $\nu(\ell(\lambda)) \geq r+1$ and $B(q)$ accounts for partitions $\mu$ satisfying $\ell(\mu)= k-1, \ \nu(\ell(\mu)) \geq j+1$ and gap between successive parts of $\mu$ being less than $r+1$. Hence the conjugates $\mu'$ of these latter partitions have the properties that $\#(\mu') = k-1, \ s(\mu') \geq j+1$, and by \eqref{freq form conj}, that each integer has a frequency less than $r+1$ in $\mu'$. Putting this information together with the partitions generated by $A(q)$ in \eqref{lemma inter 2}, the lemma follows.  
\end{proof}
 \begin{proof}[Theorem \ref{gfn for sigma rc mex}][]
By Lemma \ref{refine lemma}, we know that the set of partitions of $n$ with $j$ parts greater than the $r$-chain mex $k$ maps bijectively onto the set of partitions of $n$ whose largest $(r+1)$-repeating part is $j$, with $k$ being equal to one more than the number of parts greater than $j$. Hence, taking the sum of all the $r$-chain mex over those partitions of $n$ having $j$ parts greater than their $r$-chain mex, we obtain
\begin{equation}\label{basic block}
\sum_{\substack{\pi \in \mathcal{P}(n) \\ j \ \text{parts} > \textup{rc-mex} (\pi)}} \textup{rc-mex} (\pi) = \sum_{\substack{\pi \in \mathcal{P}(n) \\ \ell_{r+1}(\pi) =  j}} (1 + \#_{> j}(\pi)),
\end{equation}
where $\#_{> j}(\pi)$ and $\ell_{s}(\pi)$ denote respectively, the number of parts of $\pi$ that are greater than $j$ and the largest $s$-repeating part of $\pi$. Next suppose that $Q_s^{j}(n, m)$ stands for the number of partitions of $n$ whose largest $s$-repeating part is $j$, with $m$ parts being greater than $j$. Then the bivariate generating function for partitions whose largest $(r+1)$-repeating part is $j$, with $w$ keeping track of the number of parts greater than $j$ is given by
\begin{align}
\sum_{n=0}^{\infty} \sum_{m=0}^{\infty} Q_{r+1}^{j}(n, m) w^m q^n &= \frac{1}{(q; q)_{j-1}} \times \frac{q^{j(r+1)}}{1 - q^j} \times \prod_{u=j+1}^{\infty}\left(1 + wq^{u} + \cdots + w^r q^{ru}\right) \nonumber \\
&= \frac{q^{j(r+1)}}{(q; q)_j} \times \prod_{u=j+1}^{\infty} \frac{1 - (wq^u)^{r+1}}{1 - wq^u} \nonumber \\
&= \frac{q^{j(r+1)}((wq^{j+1})^{r+1}; q^{r+1})_{\infty}}{(q; q)_j(wq^{j+1}; q)_{\infty}}. \label{interm1}
\end{align}
Differentiating both sides of \eqref{interm1} with respect to $w$ and then setting $w=1$ gives us
\begin{equation} \label{diff wrt w basic}
\sum_{n=0}^{\infty} \left( \sum_{m=0}^{\infty} m Q_{r+1}^{j}(n, m)\right) q^n = \frac{q^{j(r+1)}}{(q; q)_j} \cdot \frac{d}{dw} \left.\left( \frac{((wq^{j+1})^{r+1}; q^{r+1})_{\infty}}{(wq^{j+1}; q)_{\infty}} \right)\right\vert_{w=1}.
\end{equation}
The inner sum on the left side of \eqref{diff wrt w basic} can be written as:
\begin{align}
\sum_{m=0}^{\infty} m Q_{r+1}^{j}(n, m) &= \sum_{\substack{\pi \in \mathcal{P}(n) \\ \ell_{r+1}(\pi) = j}} \#_{>j}(\pi) \nonumber \\
&= \sum_{\substack{\pi \in \mathcal{P}(n) \\ j \ \text{parts} > \textup{rc-mex} (\pi)}} \textup{rc-mex} (\pi) \quad -  \sum_{\substack{\pi \in \mathcal{P}(n) \\ \ell_{r+1}(\pi) = j}} 1 \nonumber \\
&= \sum_{\substack{\pi \in \mathcal{P}(n) \\ j \ \text{parts} > \textup{rc-mex} (\pi)}} \textup{rc-mex} (\pi) - Q_{r+1}^{j}(n), \label{sigma rc mex basic j identity}
\end{align}
where the equality in the second line above follows from \eqref{basic block} and $Q_{r}^{j}(n)$ is as defined in the proof of Theorem \ref{r-chain mex part1}. Summing the extreme sides of \eqref{sigma rc mex basic j identity} over all possible values of $j$, namely, $j \geq 0$, we arrive at
\begin{align}
\sum_{j=0}^{\infty} \sum_{m=0}^{\infty} m Q_{r+1}^{j}(n, m) &= 
\sum_{j=0}^{\infty} \sum_{\substack{\pi \in \mathcal{P}(n) \\ j \ \text{parts} > \textup{rc-mex} (\pi)}} \textup{rc-mex} (\pi) - \sum_{j=0}^{\infty} Q_{r+1}^{j}(n) \nonumber \\
\Longrightarrow \sum_{j=0}^{\infty} \sum_{m=0}^{\infty} m Q_{r+1}^{j}(n, m) &= \sigma_{rc} \textup{mex} (n) - p(n), \label{sigma rc mex basic id gfn}
\end{align} 
since the set of partitions of $n$ can be written as a disjoint union of subsets defined according to their largest $(r+1)$-repeating part, which is a non-negative integer. Finally, multiplying by $q^n$ on both sides of \eqref{sigma rc mex basic id gfn} and summing over all $n \geq 0$, we obtain 
{\allowdisplaybreaks
\begin{align}
\sum_{j=0}^{\infty} \sum_{n=0}^{\infty} \left( \sum_{m=0}^{\infty} m Q_{r+1}^{j}(n, m)\right) q^n &= \sum_{n=0}^{\infty} \sigma_{rc} \textup{mex} (n) q^n - \sum_{n=0}^{\infty} p(n) q^n \nonumber \\
\Longrightarrow \sum_{j=0}^{\infty} \frac{q^{j(r+1)}}{(q; q)_j} \cdot \frac{d}{dw} \left.\left( \frac{((wq^{j+1})^{r+1}; q^{r+1})_{\infty}}{(wq^{j+1}; q)_{\infty}} \right)\right\vert_{w=1} &= \sum_{n=0}^{\infty} \sigma_{rc} \textup{mex} (n) q^n - \frac{1}{(q; q)_{\infty}} \nonumber \\
\Longrightarrow \sum_{n=0}^{\infty} \sigma_{rc} \textup{mex} (n) q^n = \frac{1}{(q; q)_{\infty}} &+  \sum_{j=0}^{\infty} \frac{q^{j(r+1)}}{(q; q)_j} \cdot \frac{d}{dw} \left.\left( \frac{((wq^{j+1})^{r+1}; q^{r+1})_{\infty}}{(wq^{j+1}; q)_{\infty}} \right)\right\vert_{w=1}, \label{sigma rc mex final id gfn}
\end{align}
}
where the left hand side in the second step arises by using \eqref{diff wrt w basic}. We now focus on the sum in the right side of \eqref{sigma rc mex final id gfn}, namely,
\begin{equation} \label{j sum diff}
 \sum_{j=0}^{\infty} \frac{q^{j(r+1)}}{(q; q)_j} \cdot \frac{d}{dw} \left.\left( \frac{((wq^{j+1})^{r+1}; q^{r+1})_{\infty}}{(wq^{j+1}; q)_{\infty}} \right)\right\vert_{w=1}.
\end{equation}
We start with the differentiation:
\begin{align}
\frac{d}{dw} \left( \frac{((wq^{j+1})^{r+1}; q^{r+1})_{\infty}}{(wq^{j+1}; q)_{\infty}} \right) &= ((wq^{j+1})^{r+1}; q^{r+1})_{\infty} \frac{d}{dw} \frac{1}{(wq^{j+1}; q)_{\infty}} + \frac{1}{(wq^{j+1}; q)_{\infty}} \frac{d}{dw} ((wq^{j+1})^{r+1}; q^{r+1})_{\infty} \nonumber \\
&= \frac{((wq^{j+1})^{r+1}; q^{r+1})_{\infty}}{(wq^{j+1}; q)_{\infty}} \sum_{m=1}^{\infty} \frac{q^{j+m}}{1-wq^{j+m }} \nonumber \\
&- \frac{((wq^{j+1})^{r+1}; q^{r+1})_{\infty}}{(wq^{j+1}; q)_{\infty}} (r+1) w^r \sum_{m=1}^{\infty} \frac{q^{(r+1)(j+m)}}{1 - w^{r+1}q^{(r+1)(j+m)}}. \label{after diff}
\end{align}
Setting $w=1$ in \eqref{after diff}, we get
\begin{equation} \label{after w=1}
\frac{d}{dw} \left.\left( \frac{((wq^{j+1})^{r+1}; q^{r+1})_{\infty}}{(wq^{j+1}; q)_{\infty}} \right)\right\vert_{w=1} = \frac{(q^{(r+1)(j+1)}; q^{r+1})_{\infty}}{(q^{j+1}; q)_{\infty}} \left( \sum_{m = j+1}^{\infty} \frac{q^m}{1 - q^m} - (r+1) \sum_{m = j+1}^{\infty}\frac{q^{m(r+1)}}{1 - q^{m(r+1)}} \right).
\end{equation}
Putting \eqref{after w=1} in \eqref{j sum diff}, we have
\begin{align}
 &\sum_{j=0}^{\infty} \frac{q^{j(r+1)}}{(q; q)_j} \cdot \frac{d}{dw} \left.\left( \frac{((wq^{j+1})^{r+1}; q^{r+1})_{\infty}}{(wq^{j+1}; q)_{\infty}} \right)\right\vert_{w=1} \nonumber \\ 
&= \sum_{j=0}^{\infty} \frac{q^{j(r+1)}}{(q; q)_j}  \cdot \frac{(q^{(r+1)(j+1)}; q^{r+1})_{\infty}}{(q^{j+1}; q)_{\infty}} \left( \sum_{m = j+1}^{\infty} \frac{q^m}{1 - q^m} - (r+1) \sum_{m = j+1}^{\infty}\frac{q^{m(r+1)}}{1 - q^{m(r+1)}} \right) \nonumber \\
&= \frac{1}{(q; q)_{\infty}}\sum_{j=0}^{\infty} q^{j(r+1)} (q^{(r+1)(j+1)}; q^{r+1})_{\infty} \left( \sum_{m = j+1}^{\infty} \frac{q^m}{1 - q^m} - (r+1) \sum_{m = j+1}^{\infty}\frac{q^{m(r+1)}}{1 - q^{m(r+1)}} \right) \nonumber \\
&= \frac{1}{(q; q)_{\infty}} \sum_{m=1}^{\infty} \left( \frac{q^m}{1 - q^m} - (r+1) \frac{q^{m(r+1)}}{1 - q^{m(r+1)}}\right) \sum_{j=0}^{m-1} q^{j(r+1)} (q^{(r+1)(j+1)}; q^{r+1})_{\infty} \nonumber \\
&= \frac{(q^{r+1}; q^{r+1})_{\infty}}{(q; q)_{\infty}} \sum_{m=1}^{\infty} \left(\frac{q^m}{1-q^m} - (r+1) \frac{q^{m(r+1)}}{1 - q^{m(r+1)}} \right) \sum_{j=0}^{m-1} \frac{q^{j(r+1)}}{(q^{r+1}; q^{r+1})_j}. \label{inter sum simpl} 
\end{align}
Consider the inner sum in \eqref{inter sum simpl}:
\begin{equation*}
\sum_{j=0}^{m-1} \frac{q^{j(r+1)}}{(q^{r+1}; q^{r+1})_j}.
\end{equation*}
The summand here is the generating function for partitions into multiples of $r+1$ with largest part $j(r+1)$. Summing over $j$ from $0$ to $m-1$ will thus give us the generating function for partitions into multiples of $r+1$ with the largest part being at most $(m-1)(r+1)$. But this latter generating function is simply $\displaystyle\frac{1}{(q^{r+1}; q^{r+1})_{m-1}}$. Substituting this for the inner sum in \eqref{inter sum simpl} yields
\begin{align}
 &\sum_{j=0}^{\infty} \frac{q^{j(r+1)}}{(q; q)_j} \cdot \frac{d}{dw} \left.\left( \frac{((wq^{j+1})^{r+1}; q^{r+1})_{\infty}}{(wq^{j+1}; q)_{\infty}} \right)\right\vert_{w=1} \nonumber \\ 
&= \frac{(q^{r+1}; q^{r+1})_{\infty}}{(q; q)_{\infty}} \sum_{m=1}^{\infty} \left( \frac{q^m}{1-q^m} \frac{1}{(q^{r+1}; q^{r+1})_{m-1}} - (r+1) \frac{q^{m(r+1)}}{(q^{r+1}; q^{r+1})_m} \right). \label{final stretch 2 sums}
\end{align}
Now we analyze the sums corresponding to each summand in \eqref{final stretch 2 sums} separately. By an application of the $q$-binomial theorem \eqref{q-binomial}, the rightmost sum in \eqref{final stretch 2 sums} can be disposed of as follows:
\begin{align}
\sum_{m=1}^{\infty} (r+1) \frac{q^{m(r+1)}}{(q^{r+1}; q^{r+1})_m} &= (r+1) \left(-1 + \sum_{m=0}^{\infty} \frac{(q^{r+1})^m}{(q^{r+1}; q^{r+1})_m}\right) \nonumber \\
&= (r+1) \left( -1 + \frac{1}{(q^{r+1}; q^{r+1})_{\infty}}\right). \label{sec sum eval}
\end{align}
Coming to the other sum in \eqref{final stretch 2 sums}, we first express $\frac{q^m}{1 - q^m}$ as a sum of $r+1$ fractions by `going' modulo $r+1$, that is,
\begin{equation} \label{decomp}
\frac{q^m}{1 - q^m} = \frac{q^m}{1 - q^{m(r+1)}} + \frac{q^{2m}}{1 - q^{m(r+1)}} + \cdots + \frac{q^{(r+1)m}}{1 - q^{m(r+1)}}.
\end{equation}   
Substituting \eqref{decomp} back into the first summand in \eqref{final stretch 2 sums}, we get
\begin{equation}\label{first sum eval inter}
\sum_{m=1}^{\infty} \frac{q^m}{1-q^m} \cdot \frac{1}{(q^{r+1}; q^{r+1})_{m-1}} 
= \sum_{m=1}^{\infty}\left( \frac{q^m}{(q^{r+1}; q^{r+1})_{m}} + \frac{q^{2m}}{(q^{r+1}; q^{r+1})_{m}} + \cdots + \frac{q^{(r+1)m}}{(q^{r+1}; q^{r+1})_{m}}\right). 
\end{equation}
As with the computation done in \eqref{sec sum eval}, we invoke \eqref{q-binomial} for each of the $(r+1)$ summands in \eqref{first sum eval inter} and obtain
\begin{align}
\sum_{m=1}^{\infty} \frac{q^m}{1-q^m} \cdot \frac{1}{(q^{r+1}; q^{r+1})_{m-1}} &= \left( -1 + \frac{1}{(q; q^{r+1})_{\infty}}\right) +  \left( -1 + \frac{1}{(q^2; q^{r+1})_{\infty}}\right) + \cdots +  \left( -1 + \frac{1}{(q^{r+1}; q^{r+1})_{\infty}}\right) \nonumber \\
&= -(r+1) + \sum_{m=1}^{r+1} \frac{1}{(q^m; q^{r+1})_{\infty}}. \label{first sum eval}
\end{align} 
Finally, put \eqref{first sum eval} and \eqref{sec sum eval} into \eqref{final stretch 2 sums} and see that
\begin{align}
&\sum_{j=0}^{\infty} \frac{q^{j(r+1)}}{(q; q)_j} \cdot \frac{d}{dw} \left.\left( \frac{((wq^{j+1})^{r+1}; q^{r+1})_{\infty}}{(wq^{j+1}; q)_{\infty}} \right)\right\vert_{w=1} \nonumber \\
&= \frac{(q^{r+1}; q^{r+1})_{\infty}}{(q; q)_{\infty}} \left( -(r+1) + \sum_{m=1}^{r+1} \frac{1}{(q^m; q^{r+1})_{\infty}} -(r+1) \left( \frac{1}{(q^{r+1}; q^{r+1})_{\infty}} - 1\right) \right) \nonumber \\
 &= \frac{(q^{r+1}; q^{r+1})_{\infty}}{(q; q)_{\infty}} \left( \sum_{m=1}^{r} \frac{1}{(q^m; q^{r+1})_{\infty}} - \frac{r}{(q^{r+1}; q^{r+1})_{\infty}} \right). \label{pre-final}
\end{align}
Substituting \eqref{pre-final} for \eqref{j sum diff} in Equation \eqref{sigma rc mex final id gfn} finally gives us
\begin{align}
\sum_{n=0}^{\infty} \sigma_{rc} \textup{mex} (n) q^n &= \frac{1}{(q; q)_{\infty}} + \frac{(q^{r+1}; q^{r+1})_{\infty}}{(q; q)_{\infty}} \sum_{m=1}^{r} \frac{1}{(q^m; q^{r+1})_{\infty}} - \frac{r}{(q; q)_{\infty}} \nonumber \\
&= - \frac{(r-1)}{(q; q)_{\infty}} + \frac{(q^{r+1}; q^{r+1})_{\infty}}{(q; q)_{\infty}} \sum_{m=1}^{r} \frac{1}{(q^m; q^{r+1})_{\infty}} \nonumber,
\end{align}
and the proof is complete.
\end{proof}
\begin{proof}[Corollary \ref{ID FOR SIGMA RC MEX}][]
Consider a typical term for $1 \leq m \leq r$ in the sum on the right side of Theorem \ref{gfn for sigma rc mex}: 
\begin{equation} \label{InTeR}
\frac{(q^{r+1}; q^{r+1})_{\infty}}{(q; q)_{\infty}} \times \frac{1}{(q^m; q^{r+1})_{\infty}} = \frac{1}{(q^m; q^{r+1})_{\infty}^2} \prod_{\substack{n=1 \\ n \neq (u+1)(r+1), \ m + v(r+1): \ u, \ v \geq 0}}^{\infty} \frac{1}{1 - q^n}.
\end{equation}
So, the right side of \eqref{InTeR} is the generating function for partitions where no multiple of $(r+1)$ is allowed as a part, and integers congruent to $m$ modulo $(r+1)$ appear in two colors. The corollary is now immedaite.
\end{proof}

%

\section{More avenues for exploration}
Several questions present themselves for further enquiry. George Beck \cite{OEIS} conjectured and Andrews \cite{mathstu} proved a companion identity to Euler's classical result that involves the total number of parts in the two kinds of partitions. Let $a(n)$ (respectively, $b(n)$) denote the number of parts in all partitions of $n$ into odd parts (respectively, into distinct parts). Then $a(n) - b(n)$ equals the number of partitions of $n$ in which exactly one part repeats. There have been identities of a similar flavor accompanying Glaisher's identity and Franklin's identity, proved by Yang \cite{beck-type glaisher} and Ballantine and Welch \cite{beck-type franklin} respectively. This raises the question:
\begin{problem}
Find a Beck type companion identity to Proposition \ref{r-chain mex part1}, which is an analogue of Franklin's identity.
\end{problem}
We proved Theorem \ref{3-way F-analogue} by separately showing that the number of partitions of $n$ whose largest $r$-repeating is $j$ equals the other two counts of partitions. It would be worthwhile to get a direct proof between them in the spirit of Proposition \ref{r Euler mexversion}.
\begin{problem}
Find a direct proof of the identity: The number of partitions of $n$ with $j$ parts greater than $(r-1)$-chain mex equals the number of partitions of $n$ with $j$ multiples of $r$. It would also be fascinating to find a bijective proof of this identity or to begin with, the special case $r=2$, that is, Proposition \ref{r Euler mexversion}.
\end{problem}
Observe that Lemma \ref{refine lemma} is a partial refinement of Theorem \ref{3-way F-analogue}.
\begin{problem}
Find a complete refinement of Theorem \ref{3-way F-analogue}, i.e., the two stated quantities in Lemma \ref{refine lemma} also equal the number of partitions of $n$ with $j$ multiples of $r$ plus some condition depending on $k$.
\end{problem}    
Apart from a standard $q$-series proof for the generating function of $\sigma \textup{mex}(n)$, Andrews and Newman \cite[Second Proof of Theorem 1.1]{andrews-newmanI} also gave a second proof involving a combinatorial identity as a key ingredient.
\begin{problem}
It would be interesting to get a combinatorial flavored proof of Theorem \ref{gfn for sigma rc mex}.
\end{problem}  
Ballantine and Merca \cite{ballantine-merca} gave a bijective proof of \eqref{sigma mex comb id}, i.e., the identity $\sigma \textup{mex}(n) = D_2(n)$.
\begin{problem}
It is highly desirable to have a bijective proof of the identity in Corollary \ref{ID FOR SIGMA RC MEX} for $\sigma_{rc} \textup{mex}(n)$.
\end{problem}
Grabner and Knopfmacher \cite{grabner} found a Hardy-Ramanujan-Rademacher type infinite series representation for the function $\sigma \textup{mex}(n)$. We naturally ask
\begin{problem}
Find an asymptotic formula and a Hardy-Ramanujan-Rademacher type series for $\sigma_{rc} \textup{mex}(n)$.
\end{problem}
Andrews and Newman \cite{andrews-newmanI} also gave a parity characterization for $\sigma \textup{mex} (n)$, namely, that it is odd if and only if $n$ is twice a pentagonal number. 
\begin{problem}
Can we find parity results for $\sigma_{rc} \textup{mex} (n)$, at least for specific values of $r$? One may start by looking at the case of an odd $r$, since the generating function modulo $2$ of $\sigma_{rc} \textup{mex} (n)$ will be free of the term $(r-1)/(q; q)_{\infty}$.
\end{problem}
Chern \cite{Chern} introduced the maximal excludant (maex) of a partition $\pi$ to be the largest integer less than the largest part that is missing in $\pi$. In analogy with mex, we define the \emph{$t$-chain maex} of a partition $\pi$ to be the largest integer $k$ less than the largest part such that $k, \ k-1, \ \dots, \ k-t+1$ do not occur as parts in $\pi$. We call a partition $\lambda$ to be \emph{gap-free} if $s(\lambda) = 1$ and every integer between the smallest and largest parts of $\lambda$ occurs as a part. By a similar idea to that used in Theorem \ref{3-way F-analogue}, one can prove that
\begin{theorem}\label{chain maex}
The number of non-gapfree partitions of $n$ with $j$ parts greater than $(r-1)$-chain maex equals the number of partitions whose smallest $r$-repeating part is $j$.
\end{theorem}
\begin{problem}
Is there a `third side' for Theorem \ref{chain maex} similar to mex as in Theorem \ref{3-way F-analogue} which corresponds to `$j$ multiples of $r$'? It would also be worthwhile to explore  $\sigma_{rc}\textup{maex}(n)$, the sum of excludants function corresponding to $r$-chain maex.
\end{problem}

\textbf{Acknowledgements.}
The first author wants to thank the Department of Mathematics, Pt. Chiranji Lal Sharma Government College, Karnal for a conducive research environment. The second author is a SERB National Post Doctoral Fellow (NPDF) supported by the fellowship PDF/2021/001090 and would like to thank SERB for the same. The last author wants to thank SERB for the Start-Up Research Grant SRG/2020/000144.


\begin{thebibliography}{00}

\bibitem{alladi refine}
K.~Alladi, \emph{Euler's partition theorem and refinements without appeal to infinite products}, pp. 9--23 In Pillwein, V., Schneider, C. (eds) Algorithmic Combinatorics: Enumerative Combinatorics, Special Functions and Computer Algebra, Texts \& Monographs in Symbolic Computation, Springer, Cham, 2020.



\bibitem{mathstu}
G.~E.~Andrews, \emph{Euler's partition identity and two problems of George Beck}, Math. Student, ~\textbf{86} (2017), no. 1-2, 115--119.

\bibitem{gea1998}
G.~E.~Andrews, The theory of partitions, Addison-Wesley Pub. Co., NY, 300 pp. (1976). Reissued, Cambridge University Press, New York, 1998.

\bibitem{universal even}
G.~E.~Andrews and M.~Merca, \emph{On the number of even parts in all partitions of $n$ into distinct parts}, Ann. Comb., ~\textbf{24} (2020), 47--54.

\bibitem{andrews-newmanI}
G.~E.~Andrews and D.~Newman, \emph{Partitions and the minimal excludant}, Ann. Comb., ~\textbf{23} (2019), No. 2, 249--254.

\bibitem{andrews-newmanII}
G.~E.~Andrews and D.~Newman, \emph{The minimal excludant in integer partitions}, J. Integer Seq., ~\textbf{23} (2020), Article 20.2.3.

\bibitem{least r gap}
C.~Ballantine and M.~Merca, \emph{Bisected theta series, least $r$-gaps in partitions, and polygonal numbers}, Ramanujan J., ~\textbf{52} (2020), 433--444.

\bibitem{ballantine-merca}
C.~Ballantine and M.~Merca, \emph{Combinatorial proof of the minimal excludant theorem}, Int. J. Number Theory, ~\textbf{17}, No. 08 (2021), 1765--1779.

\bibitem{beck-type franklin}
C.~Ballantine and A.~Welch, \emph{Beck-type companion identities for Franklin's identity via a modular refinement}, Discrete Math., ~\textbf{344} (2021), no. 8, Paper No. 112480, 11 pp.

%

\bibitem{OEIS}
G.~Beck, The On-Line Encyclopedia of Integer Sequences, Sequence A090867, https://oeis.org.


\bibitem{Chern} S.~Chern,  \emph{Partitions and the maximal excludant},  Electron.  J.  Combin.  {\bf 28(3)} (2021),  Article P3.13

\bibitem{dhar et al}
A.~Dhar, A.~Mukhopadhyay and R.~Sarma, \emph{New relations of the mex with other partition statistics}, arXiv:2201.05997.

\bibitem{fine Euler ref}
N.~J.~Fine, Basic hypergeometric series and applications, \emph{Math. Surveys and Monographs}, ~\textbf{27} Amer. Math. Soc., Providence (1988).

\bibitem{franklin}
F.~Franklin, \emph{On partitions}, Johns Hopkins Univ. Cir., ~\textbf{2} (1883), p. 72.

\bibitem{gasper}
G.~Gasper and M.~Rahman, Basic Hypergeometric series, Second Edition, \emph{Encyclopedia of Mathematics and its applications}, 2004. 

\bibitem{grabner}
 P.J.~Grabner and A.~Knopfmacher, \emph{Analysis of some partition statistics}, Ramanujan J., ~\textbf{12} (2006), 439--454.

\bibitem{HSS}
B.~Hopkins, J.~A.~Sellers and D.~Stanton, \emph{Dyson's crank and the mex of integers partitions},  J. Comb. Theory, Ser. A, ~\textbf{185} (2022), 105523.

\bibitem{HSYee}
B.~Hopkins, J.~A.~Sellers and A.~J.~Yee, \emph{Combinatorial  perspectives on the crank and mex partition statistics}, Electron.  J.  Combin.  {\bf 29(2)} (2022),  Article P2.11


\bibitem{mex distinct}
P.~S.~Kaur, S.~C.~Bhoria, P.~Eyyunni and B.~Maji, \emph{Minimal excludant over partitions into distinct parts}, Int. J. Number Theory, ~\textbf{18}, No. 9 (2022), 2015--2028.

\bibitem{konan}
I.~Konan, \emph{A bijective proof and generalization of the non-negative crank-odd mex inequality}, arXiv:2203.04267.


\bibitem{sills}
A.~Sills, An invitation to the Rogers-Ramanujan identities, CRC Press, 2018.

\bibitem{sellers-da silva}
R.~da Silva and J.~A.~Sellers, \emph{Parity considerations for mex-related partition functions of Andrews and Newman},  J. Integer Seq., ~\textbf{23},  no. 5 (2020),
Article 20.5.7.

\bibitem{sylvester}
J.~J.~Sylvester, \emph{A constructive theory of partitions arranged in three Acts, an Interact and an Exodion}, Amer. J. Math., ~\textbf{5} (1882), 251--330.

\bibitem{beck-type glaisher}
Jane Y.~X.~Yang, \emph{Combinatorial proofs and generalizations of conjectures related to Euler's partition theorem}, European J. Combin., ~\textbf{76} (2019), 62--72.







\end{thebibliography}
\end{document}